\documentclass[12pt]{article}
\usepackage[cp850]{inputenc}
\usepackage[T1]{fontenc}\usepackage{amsfonts}
\usepackage{graphicx}
\usepackage{amsmath}
\def\<{\langle}\def\E{\mathbb{E}}\def\R{\mathbb{R}}
\def\>{\rangle}
\title{ Dirichlet random walks}
\author{G\'erard Letac\thanks{Universit\'e Paul Sabatier, 118 route de Narbonne, 31062 Toulouse, France}, Mauro Piccioni \thanks{Dipartimento di Matematica, Sapienza Universit\`{a} di Roma, 00185 Rome, Italia}}

\begin{document}\maketitle

\begin{abstract} This article provides tools for the study of the  Dirichlet random walk in $\mathbb{R}^d$. By this we mean the random variable $W=X_1\Theta_1+\cdots+X_n\Theta_n$ where $X=(X_1,\ldots,X_n) \sim \mathcal{D}(q_1,\ldots,q_n)$ is Dirichlet distributed and where  $\Theta_1,\ldots \Theta_n$ are iid, uniformly distributed on the unit sphere of $\mathbb{R}^d$ and independent of $X.$ In particular we compute explicitely in a number of cases the distribution of $W.$ Some of our results appear already  in the literature, in particular in the papers by G\'erard Le Ca\"{e}r (2010, 2011). 
In these cases, our proofs are much simpler from the original ones, since we use a kind of Stieltjes transform of $W$ instead of the Laplace transform: as a consequence  the hypergeometric functions replace the Bessel functions. 
A crucial ingredient is  a particular case of the classical and non trivial identity, true for $0\leq u\leq 1/2$:$$_2F_1(2a,2b;a+b+\frac{1}{2};u)=
\ _2F_1(a,b;a+b+\frac{1}{2};4u-4u^2).$$ 
We extend these results to a study of the limits of the Dirichlet random walks when the number of added terms goes to infinity, interpreting the results in terms of an integral by a Dirichlet process.  We introduce the ideas of Dirichlet semigroups and of Dirichlet infinite divisibility and characterize these infinite divisible distributions in the sense of Dirichlet when they are concentrated on the unit ball of $\R^d.$

\vspace{4mm}\noindent \textsc{Keywords:} Dirichlet processes, Stieltjes transforms, random flight,  distributions in a ball, hyperuniformity, infinite divisibility in the sense of Dirichlet. 

\vspace{4mm}\noindent \textsc{AMS classification}: 60D99, 60F99.
\end{abstract}

\section{Introduction} In this paper we study the distribution of $$W_{n,d}=X_1\Theta_1+\cdots+X_n\Theta_n$$ where $X=(X_1,\ldots,X_n) \sim \mathcal{D}(q_1,\ldots,q_n)$ is Dirichlet distributed and where  $\Theta_1,\ldots \Theta_n$ are iid, uniformly distributed on the unit sphere $S_{d-1}$ of the Euclidean space $E_d$  of dimension $d$ and independent of $X.$ In other terms, we select independently and uniformly $n$  random points on $S_{d-1}$ and we take their barycenter according to the random weights $(X_1,\ldots,X_n).$

We term \textit{improperly} the random variable $W_{n,d}$ a Dirichlet  random walk because if $Y_j\sim \gamma_{q_j},$ $j=1,2,\ldots$ are independent and if $S_n=Y_1+\ldots+Y_n$ then the sequence 
$$\left(\sum_{j=1}^nY_j\Theta_j\right)_{n\geq 1}$$ is a non homogeneous random walk on $E_d$ and $W_{n,d}\sim  \frac{1}{S_n}\sum_{j=1}^nY_j\Theta_j$ according to well known properties of gamma and Dirichlet distributions. 
The  random variable $W_{n,d}$ leads to a simpler analysis than the study of $\sum_{j=1}^n\Theta_j$, started by Lord Kelvin, which is the deep subject of the recent paper by Borwein \textit{et al} (2012) in the case $d=2$ and $n=3,4,5.$

If $U$ is any random variable  concentrated in the closed unit ball $B_d$ of $E_d$ our technique for characterizing the distribution of $U$ is to consider the function  defined on the interior $\stackrel{\circ}{B}_d$ of $B_d$ by $$y\mapsto T_p(U)(y)=\mathbb{E}\left(\frac{1}{(1+\<y,U\>)^{p}}\right)$$ where $p$ is a fixed positive number. The function $T_p(U)$ can be considered as a kind of Stieltjes transform of the distribution of $U.$ It is easy to prove that the knowledge of the function $T_p(U)$ gives the knowledge of the moments of the real variable  $\<y,U\>$. Since $U$ is bounded, this identifies the distribution of $U$. However, 
recovering the distribution of $U$ from $T_p(U)$ is not always easy; in particular from $T_p(U)$ and $T_q(V)$ and if $p\neq q,$ it may be puzzling to decide whether $U\sim V$ or not.

In all the cases considered here,  the distribution of $U$ is invariant by rotation, and $y\mapsto T_p(U)(y)$ depends only on $\|y\|.$ We borrow from Le Ca\"{e}r (2010) the following terminology. We say  that  $U$ is hyperspherically uniform (\textit{hyperuniform} for short) of type $k>d$ if $U$ is the law of the projection of a random variable  $\Theta$ uniform on the unit sphere   $S_{k-1}$ of $E_k$ onto a subspace of dimension $d$. The term   'hyperuniform' is justified by the fact that when  $k=d+2$ the induced distribution is uniform on the unit ball $B_d$ of $E_d$ (for $k=3$ and $d=1$ this result is due to Archimedes).


Le Ca\"{e}r's paper (2010) is a paper motivated by problems in metallurgy. It studies the cases where the Dirichlet random walk $W_{n,d}$ is hyperuniform of type $k$ for some $k>d.$ His  main result is the description of the quadruplets $(q,d,n,k)$ such that $W_{n,d}$ is hyperuniform of type $k$ when $X\sim \mathcal{D}(q,q,\ldots,q)$ (see Theorems 6 and 9 below). We refer to this paper  for the original proofs, for the bibliography concerning this problem, and for some interesting simulations: our Section 3 give new proofs of its statements.  Section 4 is devoted  to explicit expressions of the distribution of $W_{n,d}$ when $X\sim\mathcal{D}(q,\ldots,q):$  Theorem 10 deals with the case $q=d$ and gives the explicit  form of the density of $\|W_{n,d}\|^2$ as a mixing of beta distributions. Theorem 13 deals with $q=1$ and $d$ even, with results close to those of Kolesnik (2009).  Section 5 rules out the easy case $n=2.$ 
Section 6 consider the limits of $W_{n,d}$ with $(X_1,\ldots,X_n)\sim \mathcal{D}(Q/n,\ldots,Q/n)$, when $n$ goes to infinity with $Q$ and $d$ fixed, and interpret the results in terms of an integral by a Dirichlet process.  The examples of Section 6 lead to the idea of infinite divisibility in the sense of Dirichlet, introduced in Section 7. We thank A. E. Kolesnik, G. Le Ca\"{e}r,  E. Orsingher and M. Yor for interesting correspondence about the present paper.

\section {Dirichlet random walks}

\noindent \textbf{Theorem 1.} Let $(q_1,\ldots,q_n)$ be positive numbers and let $(X_1,\ldots,X_n) \sim \mathcal{D}(q_1,\ldots,q_n).$ We write $Q=q_1+\cdots+q_n.$ Let $E_d$ be the Euclidean space of dimension $d$, $\stackrel{\circ}{B}_d$ the interior of its unit ball, $S_{d-1}$ its unit sphere and let $\Theta_1,\ldots \Theta_n$ be i.i.d., each of them being uniformly distributed on $S_{d-1}.$ We consider
$$W_{n,d}=X_1\Theta_1+\cdots+X_n\Theta_n.$$ Then for $y\in \stackrel{\circ}{B}_d$ we have 
\begin{equation}\label{FAV}\mathbb{E}\left(\frac{1}{(1+\<y,W_{n,d}\>)^Q}\right)=\prod_{j=1}^{n}\ _2F_1(\frac{q_j}{2},\frac{q_j+1}{2}; \frac{d}{2}; \|y\|^2).\end{equation} 
\vspace{4mm}
\noindent
In order to give a proof we need a lemma. In the sequel $(a)_k$ is the familiar Pochhammer symbol $(a)_0=1$ and $(a)_k=a(a+1)\ldots(a+k-1)=\Gamma(a+k)/\Gamma(a)$ for $k>0.$ 

\vspace{4mm}
\noindent \textbf{Lemma 2.} Let $\Theta$ be uniform on $S_{d-1}$ and $y\in \stackrel{\circ}{B}_d$. For $p>0$ and $k$ integer we have \begin{equation}\label{MM}\mathbb{E}\left(|\<y,\Theta\>|^p\right)= \|y\|^p\frac{\Gamma(\frac{p+1}{2})\Gamma(\frac{d}{2})}{\Gamma(\frac{1}{2})\Gamma(\frac{p+d}{2})},\ \ \  \mathbb{E}\left(\<y,\Theta\>^{2k}\right)= \|y\|^{2k}(\frac{1}{2})_k/(\frac{d}{2})_k\end{equation}

\begin{equation}\label{GST}\mathbb{E}\left(\frac{1}{(1+\<y,\Theta\>)^p}\right)=\ _2F_1(\frac{p}{2},\frac{p+1}{2}; \frac{d}{2}; \|y\|^2)\end{equation}

\vspace{4mm}
\noindent \textbf{Proof of Lemma 2.} Consider the Gaussian random variable $Z=(Z_1,\ldots,Z_d)$ valued in $\mathbb{R}^d$ with distribution $N(0,I_d).$ Then $Z/\|Z\|$ is uniformly distributed on $S_{d-1}$ and $\<y,\Theta\>\sim \|y\|Z_1/\|Z\|.$ Furthermore $\<y,\Theta\>$ is a symmetric rv and the knowledge of its distribution is given by the knowledge of the distribution of its square. However
$$\frac{1}{\|y\|^2}\<y,\Theta\>^2\sim \frac{Z_1^2}{Z_1^2+\ldots+Z_d^2}\sim \beta(\frac{1}{2},\frac{d-1}{2})=\frac{1}{B(\frac{1}{2},\frac{d-1}{2})}z^{\frac{1}{2}-1}(1-z)^{\frac{d-1}{2}-1}\textbf{1}_{(0,1)}(z)dz$$ since the $Z_j^2$ are independent chi square rv with one degree of freedom. Formula (\ref{MM}) is now clear. 
Now we prove (\ref{GST}) by using the second in (\ref{MM}):
\begin{eqnarray*}
\mathbb{E}\left(\frac{1}{(1+\<y,\Theta\>)^p}\right)&=&\sum_{i=0}^{\infty}\frac{(p)_i}{i!}(-\|y\|)^i\mathbb{E}\left((\frac{1}{\|y\|}\<y,\Theta\>)^{i}\right)\\&=&\sum_{j=0}^{\infty}\frac{(p)_{2j}}{(2j)!}\|y\|^{2j}\frac{(\frac{1}{2})_j}{(\frac{d}{2})_j}=\sum_{j=0}^{\infty}\frac{(\frac{p}{2})_{j}(\frac{p+1}{2})_{j}}{(\frac{d}{2})_j j!}\|y\|^{2j}\\&=&\ _2F_1(\frac{p}{2},\frac{p+1}{2}; \frac{d}{2}; \|y\|^2)\ \square
\end{eqnarray*}

\vspace{4mm}
\noindent \textbf{Proof of Theorem 1.} We apply the following principle (see Chamayou and Letac (1994)) : if $X=(X_1,\ldots,X_n) \sim \mathcal{D}(q_1,\ldots,q_n)$ with $Q=q_1+\ldots+q_n,$ and if $f=(f_1,\ldots,f_n) $ is such that $f_i>0$ for all $i$, then
\begin{equation}\label{DM}\mathbb{E}\left(\frac{1}{\<f,X\>^Q}\right)=\frac{1}{\prod_{i=1}^nf_i^{q_i}}.\end{equation} Conditioning by $(\Theta_1,\ldots,\Theta_n)$ and applying (\ref{DM}) to $f_i=1+\<y,\Theta_i\>$ we get from (\ref{DM}) and  (\ref{GST})
\begin{eqnarray*}\mathbb{E}\left(\frac{1}{(1+\<y,W\>)^Q}\right)&=&\mathbb{E}\left(\mathbb{E}(\frac{1}{(1+\<y,W\>)^Q}|(\Theta_1,\ldots,\Theta_n))\right)\\&=&\mathbb{E}\left(\frac{1}{\prod_{i=1}^n(1+\<y,\Theta_i\>)^{q_i}}\right)\\&=&\prod_{j=1}^{n}\ _2F_1(\frac{q_j}{2},\frac{q_j+1}{2}; \frac{d}{2}; \|y\|^2).\ \ \square
\end{eqnarray*}

\vspace{4mm} \noindent In the sequel 
 we need frequently the following simple identities relating hypergeometric functions with the function 
\begin{equation}\label{FA}G(z)=\frac{2}{z}(1-\sqrt{1-z})=\frac{2}{1+\sqrt{1-z}}, |z|<1\end{equation}

\vspace{4mm} \noindent \textbf{Lemma 3.}  
We have
\begin{eqnarray}\label{HYP}\ _2F_1(\frac{c}{2}, \frac{c+1}{2}; c+1; z)&=&G(z)^c\ \mathrm{if}\ c>0,\\
\label{HYP2}\ _2F_1(\frac{c}{2}, \frac{c+1}{2}; c; z)&=&\frac{1}{\sqrt{1-z}}G(z)^{c-1}\ \mathrm{if}\ c>1,\\
\label{HYP3}\ \sum_{n=1}^{\infty}\frac{(\frac{1}{2})_n}{n!}\frac{z^n}{2n}&=&\log G(z).\end{eqnarray}

\vspace{4mm} \noindent \textbf{Proof.}
Formula (\ref{HYP}) is a consequence of $_2F_1(\frac{c}{2}, \frac{c+1}{2}; c+1; 4u-4u^2)=\ _1F_0(c;u)$ which is the formula of the abstract in a particular case. 
Formula (\ref{HYP2}) is a consequence of the Euler identity 
$$\ _2F_1(p,q;r;z)=(1-z)^{r-p-q}\ _2F_1(r-p,r-q;r;z)$$
applied to $p=c/2,$ $q=(c+1)/2$ and $r=c$,  then using (\ref{HYP}). 
If $f$ is the left hand side of (\ref{HYP3}) it is easily seen that $$f'(z)=\frac{1}{2z}\left(\frac{1}{\sqrt{1-z}}-1\right),\ f(z)=\int_0^zf'(t)dt=\int_{\sqrt{1-z}}^1\frac{ds}{1+s}=\log G(z).$$
 by the change of variable $t=1-s^2.$ $\square$

\vspace{4mm} The next proposition shows that under certain circumstances, the computation of the moments of $\|W_{n,d}\|^2$ is easy from the knowledge of $T_a(W_{n,d})(y).$

\vspace{4mm}
\noindent \textbf{Proposition 4.} Let $a,b>0.$ Let $\Theta$ uniform on the unit sphere $S_{d-1}$. Assume that $\Theta$ is independent of the random variable  $R\in[0,1],$
and denote $W=R\Theta$. Let $a>0$ and suppose that $T_a(W)$ has the form of the left hand side below for some $b>0.$ Then the even moments of $R$ have the form of the right hand side.   \begin{eqnarray}\label{FD}T_a(W)(y)=&\frac{1}{(1-\|y\|^2)^b}\ &\Rightarrow\  \mathbb{E}(R^{2k})=\frac{(b)_k(\frac{d}{2})_k}{(\frac{a}{2})_k(\frac{a+1}{2})_k}\\
\label{FD1}T_a(W)(y)=&G(\|y\|^2)^{b}\ &\Rightarrow\  \mathbb{E}(R^{2k})=\frac{(b)_{2k}(\frac{d}{2})_k}{(a)_{2k}(b+1)_k}\\
\label{FD2}T_a(W)(y)=&\frac{1}{\sqrt{1-\|y\|^2}}G(\|y\|^2)^{b-1}\ &\Rightarrow\  \mathbb{E}(R^{2k})=\frac{(b)_{2k}(\frac{d}{2})_k}{(a)_{2k}(b)_k},.\end{eqnarray}
\vspace{4mm}
\noindent \textbf{Proof.} We first expand $\mathbb{E}\left((1+\<y,W\>)^{-a}\right)$ in powers of  $\|y\|^{2}$. Using the fact that the odd moments of $\<y,W\>$ are zero and using (\ref{MM}) we get
\begin{eqnarray}\nonumber \mathbb{E}\left(\frac{1}{(1+\<y,W\>)^a}\right)&=&\sum_{k=0}^{\infty}\frac{(a)_{2k}}{(2k)!}\mathbb{E}(\<y,W\>^{2k})\\\nonumber&=&\sum_{k=0}^{\infty}\frac{(a)_{2k}}{(2k)!}\mathbb{E}(\<y,\Theta\>^{2k})\mathbb{E}(R^{2k})\\
&=&\label{EXP}\sum_{k=0}^{\infty}\frac{(a)_{2k}}{(2k)!}\frac{(\frac{1}{2})_k}{(\frac{d}{2})_k}\mathbb{E}(R^{2k})\|y\|^{2k}\end{eqnarray}
Observe also that$$(\frac{1}{2})_k=\frac{(2k)!}{2^{2k}k!},\ \ \ \  \frac{(\frac{b}{2})_k(\frac{b+1}{2})_k}{(\frac{1}{2})_k}\times \frac{(2k)!}{k!}=(b)_{2k}.$$
Using the fact that the power series of  $G(z)^b$ and $G(z)^{b-1}/\sqrt{1-z}$ are known by Lemma 3, we get the results of Proposition 4. $\square$

\vspace{4mm}
\noindent
We now apply the formulas of Proposition 4 with $R^2$ beta distributed. Actually \begin{equation}R\label{BMT}^2\sim \beta(p,q)\Leftrightarrow \mathbb{E}(R^{2k})=\frac{(p)_k}{(p+q)_k}.\end{equation}
\vspace{4mm}
\noindent \textbf{Proposition 5.} The following holds.

\begin{enumerate}\item If $W_{n,d}$ is a Dirichlet random walk governed by $\mathcal{D}(d-1,\ldots,d-1)$ and if $R_{n,d}=\|W_{n,d}\|$ then \begin{equation}\label{AFD}R_{n,d}^2\sim \beta\left(\frac{d}{2},\frac{(n-1)(d-1)}{2}\right)\end{equation} The same law is obtained if one of the Dirichlet parameters is set to $d$.

 \item If $W_{n,d}$ is a Dirichlet random walk governed by $\mathcal{D}(\frac{d}{2}-1,\ldots,\frac{d}{2}-1)$ and if $R_{n,d}=\|W_{n,d}\|$ then \begin{equation}\label{AFD1}R_{n,d}^2\sim \beta\left(\frac{d}{2},(n-1)(\frac{d}{2}-1)\right)\end{equation}
The same law is obtained if one of the Dirichlet parameters is set to $\frac{d}{2}$.
\end{enumerate}
\vspace{4mm}
\noindent \textbf{Proof.} The first statement follows from formula (\ref{FAV}) and from formula (\ref{FD}) with $a=2b=n(d-1)$ and with $a-1=2(b-1)=n(d-1)$. The second statement follows from formula (\ref{FD1}) with $a=b=n(\frac{d}{2}-1)$ and from formula (\ref{FD2}) with $a=b=n(\frac{d}{2}-1)+1$.$\square$

\vspace{4mm}
\noindent
Notice that when $d, m$ and $n$ are positive integers such that $$(n-1)(d-1)=(m-1)(d-2)=t$$ we obtain from the previous result four different representations as Dirichlet random walks for the same law $\beta(\frac{d}{2}, \frac{t}{2})$ for the squared radius. For example with $d=3, n=2, m=3$ the same law $\beta(\frac{3}{2}, 1)$ for the squared radius is obtained from 
$X \Theta_1 +(1-X) \Theta_2$, $X' \Theta_1 +(1-X') \Theta_2$, $Y_1 \Theta_1 + Y_2 \Theta_2 + Y_3 \Theta_3$ and $Y'_1 \Theta_1 + Y'_2 \Theta_2 + Y'_3 \Theta_3$, where $X \sim \beta (2,2), Y \sim \beta(2,3), X' \sim \mathcal{D}(\frac{1}{2}, \frac{1}{2}, \frac{1}{2})$ and $Y'\sim \mathcal{D}(\frac{3}{2}, \frac{1}{2}, \frac{1}{2})$. This example was already noticed in Le Ca\"{e}r (2010).

\section {Hyperuniformity}

It is easy to describe hyperuniform random vectors by introducing the standard Gaussian rv
$$Z=(Z_1,\ldots,Z_k)\sim N(0,I_k).$$
Then $U=Z/\|Z\|$ is uniformly distributed on the unit sphere of $\mathbb{R}^k$ and $(Z_1/\|Z\|,\ldots,Z_d/\|Z\|) $ is 
 hyperuniform on $E$ of type $k.$ We have 
 $$\|U\|^2=\frac{Z_1^2+\ldots+Z_{d}^2}{Z_1^2+\ldots+Z_{k}^2}\sim \beta_{d/2, (k-d)/2}$$ by standard beta gamma algebra. 
 Note that for $k-d=2$ the distribution $\beta_{d/2, 1}$ is the measure induced on the squared radius by the Lebesgue measure in $B_d,$ in agreement with Archimedes' theorem. Thus we can rephrase and extend the definition in the Introduction by allowing $k>d$ to be an arbitrary real number.
There is no particular problem in allowing $k>d$ to be a real number.
 
 \vspace{4mm} \noindent \textbf{Definition.} Let  $U$ be a random variable valued in the unit ball of the Euclidean space $E_d$ of dimension $d$  such that its distribution in invariant by rotation (namely $a(U)\sim U$ for all $a\in \mathbb{O}(E_d))$. We say that $U$ is hyperuniform of type $k>d$ if $\|U\|^2\sim \beta(d/2, (k-d)/2).$
 
 \vspace{4mm} \noindent With Proposition 5 we have immediately the  first main result of  Le Ca\"{e}r (2010):  
 

 \vspace{4mm} \noindent \textbf{Theorem 6.} Let $X=(X_1,\ldots,X_n)\sim \mathcal{D}(q+m, q,\ldots,q)$ with $m=0$ or $1$. Let $(\Theta_j)_{j=1}^n$ be iid uniform on the unit sphere of the Euclidean space of dimension $d.$ Then $W_{n,d}=X_1\Theta_1+\ldots+X_n\Theta_n$ is hyperuniform of type $k>d$ in the unit ball of $E_d$ irrespectively of $m=0,1$, if
 \begin{enumerate}\item $d\geq 2,$ $n\geq 2,$ $q=d-1$ and $k=n(d-1)+1;$
 \item $d\geq 3,$ $n\geq 2,$ $q=d/2-1$ and $k=n(d-2)+2.$ \end{enumerate}

\noindent \textbf{Proof.} 
This is immediately obtained from Proposition 5.
$\square$

 \vspace{4mm} \noindent \textbf{Proposition 7.} If $U$ in the unit ball of dimension $d$ is invariant by rotation, then $U$ is hyperuniform of type $k>d$ if and only if for any $p>0$ we have 
 $$T_p(U)(y)=\ _2F_1(p/2,(p+1)/2; k/2; \|y\|^2).$$
 
 \vspace{4mm} \noindent \textbf{Proof.} From Lemma 2 and (\ref{BMT})
 \begin{eqnarray*}\mathbb{E}\left(\frac{1}{(1+\<U,y\>)^p}\right)&=&\mathbb{E}\left(\frac{1}{(1+\<\Theta,\|U\|y\>)^p}\right)\\&=&\mathbb{E}\left(\ _2F_1(p/2,(p+1)/2; d/2; \|U\|^2\|y\|^2)\right)\\&=&\ _2F_1(p/2,(p+1)/2; k/2; \|y\|^2).\ \square\end{eqnarray*}

 \vspace{4mm} \noindent \textbf{Corollary 8.} Let $X=(X_1,\ldots,X_n)\sim D(q_1,\ldots,q_n)$ and  $Q=q_1+\ldots+q_n.$ Let $\Theta_1,\ldots,\Theta_n$ iid, uniform on the unit sphere of the Euclidean space $E$ of dimension $d$ and independent of $X.$ Then 
 $W=X_1\Theta_1+\cdots+X_n\Theta_n$ is hyperuniform if type $k>d$ if and only if for $|z|<1$ 
 
 \begin{equation}\label{EARG}
 \prod_{i=1}^n\ _2F_1(q_i/2,(q_i+1)/2; d/2; z)=\ _2F_1(Q/2,(Q+1)/2; k/2; z)
 \end{equation}
 \vspace{4mm} \noindent \textbf{Proof.} This an immediate consequence of Theorem 1 and Proposition 7, with $z$ replacing $\|y\|^2.$ $\square$

 \vspace{4mm} \noindent The second main result of Le Ca\"{e}r (2010) is a partial converse of Theorem 6 and its proof includes an unexplained miracle:

\vspace{4mm} \noindent \textbf{Theorem 9.} Let $X=(X_1,\ldots,X_n)\sim \mathcal{D}(q,\ldots,q).$  Let $(\Theta_j)_{j=1}^n$ be iid uniform on the unit sphere of the Euclidean space of dimension $d.$ Then $W_{n,d}=X_1\Theta_1+\ldots+X_n\Theta_n$ is hyperuniform of type $k>d$ in the unit ball of $E$ only if 
 \begin{enumerate}
 
 \item either $d\geq 2$, $q=d-1,$ $n\geq 2$ and $k=n(d-1)+1;$
 
 \item or $d\geq 3,$ $q=\frac{d}{2}-1,$  $n\geq 2$ and $k=n(d-2)+2.$ \end{enumerate}

\vspace{4mm} \noindent \textbf{Proof.} By Proposition 7 and Corollary 8 the random variable $W_{n,d}$ is hyperuniform of order $k$ if and only if the following hold
\begin{equation}\label{EAR}
\ _2F_1(\frac{q}{2},\frac{q+1}{2}; \frac{d}{2}; z)^{n}=\ _2F_1(\frac{nq}{2},\frac{nq+1}{2}; \frac{k}{2}; z)
 \end{equation}
With the positions $$a_i=\frac{(\frac{q}{2})_i(\frac{q+1}{2})_i}{i!(\frac{d}{2})_i}=\frac{(q)_{2i}}{2^{2i}i!(\frac{d}{2})_i}, \ \  A_i=\frac{(\frac{nq}{2})_i(\frac{nq+1}{2})_i}{i!(\frac{k}{2})_i}=\frac{(nq)_{2i}}{2^{2i}i!(\frac{k}{2})_i},$$
for $i=1,2,\ldots$, one has
$$(1+a_1z+a_2z^2)^n\equiv 1+A_1z+A_2z^2\ \ \ \mathrm{mod}\  z^3.$$
Expanding the left hand side member through the multinomial formula we get the two equations
\begin{eqnarray}
\label{EAR1}A_1&=&na_1\\
\label{EAR2}A_2&=&na_2+\frac{n(n-1)}{2}a_1^2\\
\end{eqnarray}
Equality (\ref{EAR1}) implies $k(q+1)=d(nq+1).$ Equality (\ref{EAR2}) after obvious simplifications becomes
$$\frac{(nq+1)(nq+2)(nq+3)}{k(k+2)}=\frac{(q+1)(q+2)(q+3)}{d(d+2)}+(n-1)\frac{q(q+1)^2}{d^2}$$
Replacing $k$ by $d(nq+1)/(q+1)$ in this equation and simplifying again leads to 
$$\frac{(q+1)(nq+2)(nq+3)}{d(nq+1)+2q+2}=\frac{(q+2)(q+3)}{d+2}+(n-1)\frac{q(q+1)}{d}$$
Next, set $x=\frac{d}{q+1}$, obtaining (for fixed $q$) an equation in $x$
\begin{equation}\label{EMI}\frac{(nq+2)(nq+3)}{x(nq+1)+1}=\frac{(q+2)(q+3)}{x(q+1)+1}+(n-1)\frac{q}{X}\end{equation}
Now we attend to a little \textit{miracle} of the theory: the equation (\ref{EMI}) is of second degree equation in $x$, whose solutions $x=1$ and $x=2$ depend neither on $n$ nor in $q.$ From this we get that (\ref{EMI}) has the only solutions $q=d-1$ and $q=\frac{d}{2}-1$, as desired. $\square$

\vspace{4mm} \noindent \textbf{Comment.} From Theorems 6 and   9 the only Dirichlet walks $W_{n,d}$ governed by $(X_1,\ldots,X_n)\sim \mathcal{D}(q,\ldots,q)$ are uniform (we mean that they  are hyperuniform with type $k=d+2)$ if and only if 
 \begin{enumerate}
 \item $d=2,$ $n=3,$ $q=1.$ 
 \item $d=3,$ $n=2,$ $q=2.$ 
 \item $d=3,$ $n=3,$ $q=1/2.$ 
 \item $d=4,$ $n=2,$ $q=1.$
\end{enumerate}

\section{Density of the Dirichlet random walk in the  case $\mathcal{D}(q,\ldots,q)$ }  

Let us study the density of 
the squared radius of the Dirichlet random walk $W_{n,d}$ governed by $\mathcal{D}(q,\ldots,q)$, with $n \geq 2$. The cases $q=d-1$ and $q=\frac{d}{2}-1$ have been already dealt with in Theorem 6. First, we rule out the simple case $d=1$ for arbitrary $q$ (Proposition 10).  Then we present the main results of the present section for the case $q=d$ (Theorem 11) and the case $q=1$ and $d$ even (Theorem  13). The case $q=d$  is the subject of the second paper of Le Ca\"{e}r (2011).  The case $q=d=2$ is considered by Beghin  and Orsingher (2010).  The case $(q,d)=(1,6)$ is the subject of  Kolesnik (2009).

\vspace{4mm} \noindent \textsc{Case $d=1$, arbitrary $q>0.$} 

\vspace{4mm} \noindent
For this case we can deal  with a slightly more general case. Dirichlet random walk occurs when $p=1/2$ in the following proposition:

\vspace{4mm} \noindent \textbf{Proposition 10.} Let $X=(X_1,\ldots,X_n)\sim \mathcal{D}(q,\ldots,q)$ ($n$ times) and let $\epsilon_1,\ldots,\epsilon_n$ be iid random variables independent of $X$ with distribution $p\delta_1+(1-p)\delta_{-1}.$ Then the law of $Y=X_1\epsilon_1+\cdots+X_n\epsilon_n$ is 
$$(1-p)^n\delta_{-1}(dy)+p^n\delta_1(dy)+\Gamma(nq)2^{-nq-1}\sum_{k=1}^{n-1}\binom{n}{k}p^{k}(1-p)^{n-k}\frac{(1+y)^{kq-1}(1-y)^{(n-k)q-1}}{\Gamma(kq)\Gamma((n-k)q)}\textbf{1}_{(-1,1)}(y)dy.$$ 


\vspace{4mm} \noindent \textbf{Proof.} Easy by conditioning on $K=\epsilon_1+\ldots+\epsilon_n.$

  \vspace{4mm} \noindent \textsc{Case $q=d>1$.}

\vspace{4mm} \noindent With $R_{n,d}=\|W_{n,d}\|$ as usual, Theorem 1 and formula (\ref{FD}) show that for $q=d$
\begin{equation}\label{MTDQ}\mathbb{E}(R^{2k})=\frac{(\frac{nd+n}{2})_k(\frac{d}{2})_k}{(\frac{nd}{2})_k(\frac{nd+1}{2})_k}=\frac{(\frac{nd+n}{2})_k}{(\frac{nd+1}{2})_k}\times\frac{(\frac{d}{2})_k}{(\frac{nd}{2})_k}\end{equation}
 In other terms, if $Z\sim R^2$ and $X$ are independent, then $X\sim \beta( \frac{nd+1}{2},\frac{n-1}{2})$ implies 
\begin{equation}\label{MCE} Y=XZ \sim \beta(\frac{d}{2}, \frac{(n-1)d}{2}). \end{equation}
This is a consequence of (\ref{BMT}). However, finding the distribution of $Z$ from  the multiplicative convolution equation (\ref{MCE}) is not easy. Here is the result:

\vspace{4mm} \noindent \textbf{Theorem 11.} Suppose that $d>1.$ Consider  the  Dirichlet random walk $W_{n,d}$ governed by $\mathcal{D}(d,\ldots,d)$, denote  $R_{n,d}=\|W_{n,d}\|$ and let  $f$ be the density of $V=R_{n,d}^2.$ Then the Mellin transform $M$ of $f$ is given by 
\begin{equation}\label{MTLC}M(s)=\int_0^{1}f(v)v^{s}dv=C\frac{\Gamma(\frac{nd+n}{2}+s)\Gamma(\frac{d}{2}+s)}{\Gamma(\frac{nd+1}{2}+s)\Gamma(\frac{nd}{2}+s)}\end{equation} where $C$ is determined by $M(0)=1.$ Furthermore
\begin{enumerate}\item if $n=2N+1$ is odd, then $M$ is a rational function of the form
$$M(s)=C\frac{(s+Nd+\frac{d+1}{2})_N}{(\frac{d}{2}+s)_{Nd}}=\sum_{k=0}^{Nd-1}\frac{A_k}{\frac{d}{2}+k+s}$$
and $f(v)=\sum_{k=0}^{Nd-1}A_kv^{\frac{d}{2}+k-1}, $ where the sign of $A_k$ is $(-1)^k;$
\item if $n=2N$ is even, then  
$$M(s)=C\frac{\Gamma(\frac{d}{2}+s)}{\Gamma(Nd+\frac{1}{2}+s)}\times (s+Nd)_N$$
and $f$ is a mixture of beta densities:
$$f(v)dv=\sum_{k=0}^{N}r_k\beta(\frac{d}{2}+k,Nd-\frac{d-1}{2}-k)(dv)$$ with the weights $r_0,\ldots,r_N$ 
that are positive, satisfy $\sum_{k=0}^Nr_k=1$  and have the explicit value \begin{equation}\label{VCB}r_k=\frac{1}{(Nd)_N}\left(\frac{d}{2}\right)_k\left(d(N-\frac{1}{2})\right)_{N-k}\binom{N}{k},\,\, k=0,\ldots,N.\end{equation}

\end{enumerate}
\vspace{4mm} \noindent \textbf{Proof.} With the notation $X,Y,Z$ used in  (\ref{MCE}) we  can claim that 
$\mathbb{E}(X^s)\mathbb{E}(Z^s)=\mathbb{E}(Y^s)$ for $s>0$ and this implies (\ref{MTLC}). When $n=2N+1$ is odd the $A_k$ are computed by partial fraction expansion and $\frac{1}{\frac{d}{2}+k+s}$ is the Mellin transform of the function $v^{\frac{d}{2}+k-1}\textbf{1}_{(0,1)}(v).$  

When $n=2N$ the situation is more complicated. Since $(\frac{(x)_k}{k!})_{k\geq 0}$ is a basis for real polynomials, every polynomial $P(x)$ of degree $N$ can be written as $$P(x)=\sum_{k=0}^Np_k\frac{(x)_k}{k!}$$ for some coefficients $p_0,\ldots,p_N.$ This implies that we have
$
P(-j)=\sum_{k=0}^j(-1)^kp_k\binom{j}{k}$ 
 which implies $p_k=\sum_{j=0}^k(-1)^jP(-j)\binom{k}{j}.$ 
We now apply these remarks to the polynomial $P(x)= \frac{(x+A)_N}{(Nd)_N},$ with $A=d(N-\frac{1}{2})$. Taking $x=s+\frac{d}{2}$, we get
$P(s+\frac{d}{2})=\frac{(s+Nd)_N}{(Nd)_N},$ thus 
\begin{equation}\label{RL}\frac{(s+Nd)_N}{(Nd)_N}=\sum_{k=0}^Nr_k\frac{(s+\frac{d}{2})_k}{(\frac{d}{2})_k}\end{equation}
We proceed to   an explicit calculation of the $r_k$'s as follows:  
$$r_k=\frac{1}{(Nd)_N}\left(\frac{d}{2}\right)_k\sum_{j=0}^k\frac{(-1)^j}{j!}\frac{1}{(k-j)!}(A-j)_N=
\frac{1}{(Nd)_N}\left(\frac{d}{2}\right)_k(A)_{N-k}a_k$$
where 
$$a_k=\sum_{j=0}^k(-1)^j\frac{(A-j)_j}{j!}\frac{(A+N-k)_{k-j}}{(k-j)!}$$ 
For computing the $a_k$'s we use their generating function and the change of index $k'=k-j:$
\begin{eqnarray*}\sum_{k=0}^{\infty}a_kz^k&=&\sum_{k=0}^{\infty}\sum_{j=0}^k(-1)^j\frac{(A-j)_j}{j!}z^j\frac{(A+N-k)_{k-j}}{(k-j)!}z^{k-j}\\&=&\sum_{j=0}^{\infty}(-1)^j\frac{(A-j)_j}{j!}z^j\sum_{k'=0}^{\infty}\frac{(A+N-j-k')_{k'}}{k'!}z^{k'}\\&=&\sum_{j=0}^{\infty}(-1)^j\frac{(A-j)_j}{j!}z^j(1+z)^{A+N-j}\\&=&(1+z)^{A+N-1}\sum_{j=0}^{\infty}\frac{(A-j)_j}{j!}\left(\frac{-z}{1+z}\right)^j\\&=&(1+z)^{A+N-1}\left(1-\frac{z}{1+z}\right)^{A-1}=(1+z)^N
\end{eqnarray*}
which proves $a_k=\binom{N}{k}$ and (\ref{VCB}). Finally plugging $s=0$ in (\ref{RL}) we get that $\sum_{k=0}^Nr_k=1.$ 
To end up the proof, we have now from (\ref{MTLC}) and (\ref{RL})
$$M(s)=C\sum_{k=0}^Nr_k\frac{\Gamma(\frac{d}{2}+s)}{\Gamma(Nd+\frac{1}{2}+s)}\frac{(s+\frac{d}{2})_k}{(\frac{d}{2})_k}$$
where $C=\frac{\Gamma(Nd+\frac{1}{2})}{\Gamma(\frac{d}{2})}$ and we observe that 
$$C\frac{\Gamma(\frac{d}{2}+s)}{\Gamma(Nd+\frac{1}{2}+s)}\frac{(s+\frac{d}{2})_k}{(\frac{d}{2})_k}=\frac{\Gamma(Nd+\frac{1}{2})}{\Gamma(\frac{d}{2}+k)}\frac{\Gamma(\frac{d}{2}+k+s)}{\Gamma(Nd+\frac{1}{2}+s)}$$ is the Mellin transform of $\beta(\frac{d}{2}+k, Nd-\frac{d-1}{2}-k).$ $\square$

\vspace{4mm} \noindent \textsc{Case $d$ even
 and $q=1$.}

\vspace{4mm} \noindent The next theorem (Theorem 13) is inspired by Kolesnik (2009) who computes the distribution of $W$ for $d=6$ when $W$ is a Dirichlet random walk governed by $\mathcal{D}(1,1,\ldots,1)$ but when $n$ is random and Poisson distributed.  More specifically, let $Y_0,\ldots,Y_n,\ldots$ be iid rv such that $\Pr(Y_0>y)=e^{-cy}$, let $N(t)$ be a Poisson process such that $\mathbb{E}(N(t))=\lambda t$ and let  $\Theta_0,\ldots,\Theta_n,\ldots$ be iid and uniform on the unit sphere of a Euclidean space of dimension $d.$ Denote
\begin{equation}\label{POK}X(t)=\sum_{i=0}^{N(t)}Y_i\Theta_i,\ S(t)=\sum_{i=0}^{N(t)}Y_i,\ W(t)=X(t)/S(t).\end{equation}
Kolesnik computes the distribution of $X(t)$ for $d=6.$ Here we are rather interested in the distribution of $W(t)|N(t)=n-1$, which is independent of $S(t)|N(t)=n-1.$

Before stating Theorem 13, we need a presentation of the hypergeometric function $\  _2F_1(\frac{1}{2},1;D;z)$ when $D\geq 2$ is a positive integer. The basic point of Proposition 12 is that this function is actually a polynomial $\sum_{k=1}^{D-1}B_kG(z)^k$ with respect to the function $G$ defined by (\ref{FA}).

\vspace{4mm} \noindent \textbf{Proposition 12.} Let  $D$ be a positive integer $\geq 2.$  Then for $0<u<1/2$ 
\begin{equation}\label{BKO}  _2F_1(\frac{1}{2},1;D;4u-4u^2)=\frac{(D-1)!B_D(u)}{4^{D-1}(1-u)^{D-1}}=\sum_{k=1}^{D-1}\frac{B_k}{(1-u)^k}\end{equation} where $B_D(u)$ is a polynomial of degree $D-2$ defined as follows: if $A_D(u)=u^{D-1}B_D(u)$ then $A'_{D}(u)=4(1-2u)A_{D-1}(u)$, $A_D(0)=0$ and $A_1(u)=\frac{1}{1-2u}.$ In particular  $B_1=1$ for $D=2,$ $B_1=4/3$ and $B_2=-1/3$
for $D=3$ and $B_1=3/5,$ $B_2=-6/5$  and $B_3=8/5$ for $D=4.$

\vspace{4mm} \noindent \textbf{Proof.} $$\  _2F_1(\frac{1}{2},1;D;z)=(D-1)!\frac{1}{z^{D-1}}H_D(z)$$ where 
$$H_D(z)=\sum_{n=0}^{\infty}\frac{(\frac{1}{2})_n}{(n+D-1)!}z^{n+D-1}.$$ Therefore $H_D(0)=H_D'(0)=\ldots=H_D^{(D-2)}(0)=0$ and for $0\leq i<D$ we have 
$$(\frac{d}{dz})^{i}H_D(z)=H_{D-i}(z),\ \ (\frac{d}{dz})^{D-1}H_D(z)=H_1(z)=\frac{1}{\sqrt{1-z}}.$$ 
We now show by induction on $D>1$ that $u\mapsto A_D(u)= H_D(4u-4u^2)$ is a polynomial of degree $ 2D-3$ and of valency $D-1.$ Note that  $A_1(u)=H_1(4u-4u^2)=1/(1-2u).$ Since for $D\geq 2$  $A'_{D}(u)=4(1-2u)A_{D-1}(u)$ and since $A_D(0)=0$ then $A_2(u)=4u$ and the property is true for $D=2.$  

Suppose that the property is true for some $2\leq D.$ Then $A'_{D+1}(u)=4(1-2u)A_{D}(u)$ and  the induction hypothesis shows that $A_{D+1}$ is a polynomial of degree $2D-1.$  Introduce the polynomial $B_D$ such that $A_D(u)=u^{D-1}B_D(u).$ It exists from the induction hypothesis. Since $A_{D+1}(0)=0$ we have  $A_{D+1}(u)=4\int_0^u(1-2v)v^{D-1}B_{D}(v)dv$. Therefore the valency of $A_{D+1}$ is $D$ and the induction hypothesis is extended. The remainder is plain. $\square$

\vspace{4mm} \noindent \textbf{Theorem 13.} Let $W_{n,2D}$  be a Dirichlet random walk in the unit ball of $\mathbb{R}^{2D}$ governed by $\mathcal{D}(1,1,\ldots,1)$ ($n$ times) where $D$ is a positive integer $\geq 2.$ Define the sequence $\{p_i,\  i=n,n+1, \ldots,n(D-1)\}$ by 
$$\left(\sum_{k=1}^{D-1}B_kz^k\right)^n=\sum_{i=n}^{n(D-1)}p_iz^i$$ where $B_1,\ldots, B_{D-1}$ are the numbers defined by (\ref{BKO}). Then 
$$T_n(W_{n,d})(y)=\sum_{i=n}^{n(D-1)}p_i\ _2 F_1(\frac{i}{2},\frac{i+1}{2};i+1;\|y\|^2)$$
In particular, if $R_{n,d}=\|W_{n,d}\|$ we have  the moments of $R_{n,d}^2$  
\begin{equation}\label{KOD}\mathbb{E}(R_{n,d}^{2k})=\frac{(D)_k}{(n)_{2k}}\sum_{i=n}^{n(D-1)}p_i\frac{(i)_{2k}}{(i+1)_k}
\end{equation}
as well as the Mellin transform of $R_{n,d}^2$
\begin{equation}\label{KOD1}\mathbb{E}(R_{n,d}^{2s})=\frac{(n-1)!}{(D-1)!}\sum_{i=n}^{n(D-1)}ip_i\frac{(s)_D}{(s)_{i+1}}\frac{(2s)_i}{(2s)_n}
\end{equation}

\vspace{4mm} \noindent \textbf{Proof.} From Theorem 1 we  have 
$$\mathbb{E}\left(\frac{1}{(1+\<y,W_{n,d}\>)^n}\right)=\left(\  _2F_1(\frac{1}{2},1;D;\|y\|^2)\right)^n$$
From Proposition 12 we have 
$$\left(\  _2F_1(\frac{1}{2},1;D;4u-4u^2)\right)^n=\sum_{i=n}^{n(D-1)}p_i\frac{1}{(1-u)^i}.$$ 
With the notation (\ref{FA}) for $G$ and using formula (\ref{HYP}) we get that 
\begin{eqnarray*}\left(\  _2F_1(\frac{1}{2},1;D;z)\right)^n&=&\sum_{i=n}^{n(D-1)}p_i G(z)^i\\&=&\sum_{i=n}^{n(D-1)}p_i\ _2F_1(\frac{i}{2},\frac{i+1}{2};i+1;z)\\&=&\sum_{k=0}^{\infty}\frac{z^k}{k!}\left(\sum_{i=n}^{n(D-1)}p_i\frac{(\frac{i}{2})_k(\frac{i+1}{2})_k}{(i+1)_k}\right)=\sum_{k=0}^{\infty}\frac{z^k}{k!}\left(\sum_{i=n}^{n(D-1)}p_i\frac{(i)_{2k}}{2^{2k}(i+1)_k}\right)\end{eqnarray*}
Next we write write $W_{n,d}=R_{n,d}\Theta$ and use (\ref{EXP}) to come up with formula (\ref{KOD}). For proving 
(\ref{KOD1}) we rewrite (\ref{KOD}) as follows: 
\begin{eqnarray*}\mathbb{E}(R_{n,d}^{2k})&=&\sum_{i=n}^{n(D-1)}p_i\frac{\Gamma(D+k)}{\Gamma(D)\Gamma(k)}\times \frac{\Gamma(i+1)\Gamma(k)}{\Gamma(i+1+k)}\times \frac{\Gamma(i+2k)}{\Gamma(i)\Gamma(2k)}\times \frac{\Gamma(n)\Gamma(2k)}{\Gamma(n+2k)}\\&=&\frac{(n-1)!}{(D-1)!}\sum_{i=n}^{n(D-1)}ip_i\frac{(k)_D}{(k)_{i+1}}\frac{(2k)_i}{(2k)_n}\end{eqnarray*} Therefore (\ref{KOD1}) is correct when $s$ is a non negative integer. Now observe that the right hand side $H(s)$ of (\ref{KOD1}) is a rational fraction. More specifically, it is a linear combination  of the  rational fractions $\frac{(s)_D}{(s)_{i+1}}\frac{(2s)_i}{(2s)_n}$ and the difference between the degree of the denominator and that of the numerator is $n+1-D$ and does not depend on $i.$  If $n+1>D$ the rational fraction $H$ is a linear combination of partial fractions of the form $\frac{1}{s+j}=\int_0^1 v^{s+j-1}dv$, and $H$ is the Mellin transform of a polynomial $P_H$ restricted to $(0,1).$ Since $\mathbb{E}(R_{n,d}^{2k})=\int_0^1v^kP_H(v)dv$ this implies that $$R_{n,d}^2\sim P_H(v)\textbf{1}_{(0,1)}(v)dv$$ and therefore (\ref{KOD1}) holds for all $s\geq 0.$ If $n+1\leq D$ the rational fraction $H$ is the sum of a polynomial $Q$ of degree $\leq D-n-1$ and of a linear combination of partial fractions of the form $\frac{1}{s+j}=\int_0^1 v^{s+j-1}dv.$ We claim that $Q$ is the zero polynomial. If $\deg Q>0$ then  $k\mapsto \mathbb{E}(R^{2k})$ is unbounded, which is impossible since $0\leq R_{n,d}\leq 1.$ If $Q$ is a  non zero constant $c$, this implies 
that $\lim_{k\rightarrow \infty \infty}\mathbb{E}(R_{n,d}^{2k})=c$ which means that $\Pr(R_{n,d}^2=1)=c$. This fact  is impossible since $R_{n,d}=\|W_{n,d}\|=\|X_1\Theta_1+\ldots+X_n\Theta_n\|=1$ implies that all $\Theta_j$ are equal, which has zero probability  if $d>1.$ Therefore $Q=0$ and one concludes that (\ref{KOD1}) holds as in the case $n+1>D$. $\square$

\vspace{4mm} \noindent \textbf{Example.} 
\textsc{Case $q=1,$ $d=6$ and $n=2.$} We have seen in Proposition 12 that in that case, since $D=d/2=3$ we have 
$$_2F_1(\frac{1}{2},1;3;z)=\frac{4}{3}G(z)-\frac{1}{3}G^2(z)$$
which implies $$\left(\ _2F_1(\frac{1}{2},1;3;z)\right)^2=\frac{16}{9}G^2(z)-\frac{8}{9}G^3(z)+\frac{1}{9}G^4(z)$$
that is to say $p_2=\frac{16}{9},\ p_3=-\frac{8}{9},\ p_4=\frac{1}{9}.$ Careful calculations from (\ref{KOD1}) give
$$\mathbb{E}(R_{6,2}^{2s})=\frac{8}{3+s}-\frac{20}{3}\frac{1}{4+s}=\int_0^1 v^s(8v^2-\frac{20}{3}v^3)dv$$
In other terms $R_{6,2}\sim (8v^2-\frac{20}{3}v^3)\textbf{1}_{(0,1)}(v)dv.$  This is equivalent to formula (13) of Kolesnik (2009). Applying the same method as above for $n=3$ would provide formula (15) of Kolesnik.

\section {The case $\mathcal{D}(q_1,q_2)$} Since the Dirichlet distribution $\mathcal{D}(q_1,q_2)$ is nothing but the distribution of $(X,1-X)$ where $X\sim \beta(q_1,q_2),$ it is almost trivial to study directly the  Dirichlet random walks for $n=2,$ but it offers the opportunity to check  general formulas in this particular case. 

\vspace{4mm} \noindent \textbf{Proposition 14.} Let $X\sim \beta(q_1,q_2)$, $\Theta_1$ and $\Theta_2$ three independent random variables such that the $\Theta_i$ are uniformly distributed on the unit sphere $S$ of the Euclidean space of dimension $d.$ Let $R_{2,d}=\|W_{2,d}\|$ where $W_{2,d}$ is the Dirichlet random walk 
\begin{equation}\label{TWO}W_{2,d}=X\Theta_1+(1-X)\Theta_2.\end{equation} Then the Mellin transform of $H_{2,d}=1-R_{2,d}^2$ is 
$$\mathbb{E}(H_{n,d}^s)=C\frac{\Gamma(q_1+s)\Gamma(q_2+s)\Gamma(\frac{d-1}{2}+s)}{\Gamma(\frac{q_1+q_2}{2}+s)\Gamma(\frac{q_1+q_2+1}{2}+s)\Gamma(d-1+s)}$$
where $C$ is the normalizing constant  such that the right hand side is 1 when $s=0.$

\vspace{4mm} \noindent \textbf{Proof.} Setting $Z=\frac{1}{2}(1-\<\Theta_2,\Theta_1\>)$ we have 
$$1-H_{2,d}=R_{n,d}^2\sim \|X\Theta_1+(1-X)\Theta_2\|^2=1-4X(1-X)Z.$$ Since $X$ and $Z$ are independent we have $\mathbb{E}(H_{2,d}^s)=\mathbb{E}((4X(1-X))^s)\mathbb{E}(Z^s).$ Since $X\sim \beta(q_1,q_2)$ then 
$$\mathbb{E}((4X(1-X))^s)=C_1\frac{\Gamma(q_1+s)\Gamma(q_2+s)}{\Gamma(\frac{q_1+q_2}{2}+s)\Gamma(\frac{q_1+q_2+1}{2}+s)}$$with a suitable normalizing constant $C_1.$ From Lemma 2, conditioning first on $\Theta_2$, we know that $\<\Theta_2,\Theta_1\>$ is symmetric and that $\<u,\Theta_1\>^2\sim \beta(\frac{1}{2},\frac{d-1}{2})$; this implies easily that $Z\sim \beta(\frac{d-1}{2},\frac{d-1}{2}).$ From this the computation of $\mathbb{E}(Z^s)$ can be done and this leads  to the proof of the proposition. $\square$

\vspace{4mm} \noindent \textbf{Example.} If  we apply the above proposition to  $\mathcal{D}(q,q)$ and to $\mathcal{D}(q,q+1)$ we obtain in both cases the same result
$$\mathbb{E}(H_{2,d}^s)=\frac{\Gamma(d-1)\Gamma(q+\frac{1}{2})\Gamma(q+s)\Gamma(\frac{d-1}{2}+s)}{\Gamma(q+\frac{1}{2}+s)\Gamma(d-1+s)\Gamma(q)\Gamma(\frac{d-1}{2})}.$$ 
The above equality says that $H_{2,d}\sim XY$ where $X\sim \beta(\frac{d-1}{2},\frac{d-1}{2})$ and $Y\sim \beta(q,\frac{1}{2})$ are independent. 
When $q=d-1$
$$\mathbb{E}(H_{2,d}^s)=\frac{\Gamma(d-\frac{1}{2})}{\Gamma(\frac{d-1}{2})}\frac{\Gamma(\frac{d-1}{2}+s)}{\Gamma(d-\frac{1}{2}+s)}$$ and therefore $H_{2,d}\sim \beta(\frac{d-1}{2},\frac{d}{2})$ and $R_{2,d}^2\sim \beta(\frac{d}{2},\frac{d-1}{2}), $ in agreement with Proposition 5.

\vspace{4mm} \noindent \textbf{Corollary 15.} If $X$ is uniform on $(0,1)$ the Mellin transform of $H_{2,d}=1-\|W_{2,d}\|^2$ as defined by (\ref{TWO}) is as follows \begin{enumerate}

\item if $d=2$ then $\mathbb{E}(H_{2,d}^s)=1/(1+2s)$ which implies $H_{2,d}\sim \frac{1}{2h^{1/2}}\textbf{1}_{(0,1)}(h)dh$ and $R_{2,d}^2 \sim \beta(1,\frac{1}{2})$.

\item if $D \geq 2$ then $$\mathbb{E}(H_{2,2D}^s)=\frac{(2D-2)!}{(\frac{3}{2})_{D-2}}\frac{(s+\frac{3}{2})_{D-2}}{(s+1)_{2D-2}};$$
\item if $D \geq 1$ then $$\mathbb{E}(H_{2,2D+1}^s)=\frac{\Gamma(\frac{3}{2})(2D-1)!}{(D-1)!}\times \frac{\Gamma(s+1)}{\Gamma(s+\frac{3}{2})}\times \frac{1}{(s+D)_D}.$$
\end{enumerate}

\vspace{4mm} \noindent \textbf{Example.} If $d=6$ part 2 shows that \begin{eqnarray*}\mathbb{E}(H_{2,6}^s)&=&16\frac{s+\frac{3}{2}}{(s+1)(s+2)(s+3)(s+4)}\\&=&\frac{20}{3}\frac{1}{s+4}-\frac{12}{s+3}+\frac{4}{s+2}+\frac{4}{3}\frac{1}{s+1}\\&=&\int_0^1\left(\frac{20}{3}h^3-12h^2+4h+\frac{4}{3}\right)h^sdh\end{eqnarray*}
Therefore the distribution of $R_{2,6}^2=1-H_{2,6}$ is $$\left(\frac{20}{3}(1-v)^3-12(1-v)^2+4(1-v)+\frac{4}{3}\right)\textbf{1}_{(0,1)}(v)dv=(8v^2-\frac{20}{3}v^3)\textbf{1}_{(0,1)}(v)dv$$ This result of course coincides with the result of the last example of Section 4.

\section{Limits of  Dirichlet walks}
It is quite natural to ask what happens in Theorem 1 as $n\rightarrow\infty$ and $q_j=Q/n$ for $j=1,\ldots,n$, for some positive constant $Q>0$. The answer turns out to be rather surprising.

\vspace{4mm} \noindent \textbf{Theorem 16.} Let $q_j=Q/n$ for $j=1,\ldots,n$ in the definition of $W_{n,d}$. Then the following hold:
\begin{enumerate}
\item If $n\rightarrow\infty$, then for $y \in \stackrel{\circ}{B}_d$ we have $T_Q(W_{n,d})(y) \rightarrow e^{L_d(\|y\|^2)}$
where 
$$L_d(z)=\frac {1}{2}\sum_{k=1}^{\infty}\frac{(1/2)_k}{(d/2)_k}\frac{z^k}{k}.$$

\item For $y \in \stackrel{\circ}{B}_d$ and $\Theta$ uniform on the sphere $S_{d-1}$ denote by $U_1$ the first component of $\Theta.$ Then $$L_d(\|y\|^2)=-\E\{\log(1+\|y\|U_1)\}.$$

\item $W_{n,d}$ tends weakly to $W^Q_{\infty,d}$, which is characterized by
$$ T_Q(W^Q_{\infty,d})(y)=\exp\{-Q\E\{\log(1+\|y\|U_1)\}$$

\item Let $\{\Theta_i, i=1,2,\ldots\}$ is a sequence of iid rv's uniform on the unit sphere $S_{d-1}$ and $\{\Pi_i, i=1,2,\ldots\}$ be obtained from an independent iid sequence  $\{Y_i, i=1,2,\ldots\}$ of $\beta (1,Q)$ random variables, through
$\Pi_1=Y_1$ and $\Pi_n=Y_n \prod_{i=1}^{n-1} (1-Y_i)$. Then
$$W^Q_{\infty,d} \sim \sum_{i=1}^{\infty} \Pi_i \Theta_i.$$
\end{enumerate}
\vspace{4mm} \noindent \textbf{Proof.} 
Statement 1 results from Theorem 1 and the development
$$_2F_1(\frac{Q}{2n},\frac{Q}{2n}+\frac{1}{2};\frac{d}{2};z)=1+\frac {Q}{n}L_d(z)+o(\frac{1}{n})$$ from which we get 
$$\left(\ _2F_1(\frac{Q}{2n},\frac{Q}{2n}+\frac{1}{2};\frac{d}{2};\|y\|^2)\right)^n \rightarrow e^{QL_d(\|y\|^2)}$$
On the other hand, recalling that $W_{n,d}=\sum_{i=1}^n X_i^{(n)} \Theta_i$, with $X^{(n)}=(X_1^{(n)},\ldots, X_n^{(n)}) \sim \mathcal{D}(\frac{Q}{n},\ldots,\frac{Q}{n})$
one has for $y \in \stackrel{\circ}{B}_d$
$$
T_Q(W_{n,d})(y) =\E \left(\frac{1}{(\sum_{i=1}^nX_i^{(n)}(1+\<y,\Theta_i\>))^Q}\right)
$$
which, by using (\ref{DM}), can be written as 
$$
 \E\left( \frac{1}{\prod_{i=1}^n (1+\<y,\Theta_i\>)^{\frac{Q}{n}}} \right)
=\left(\E(\frac{1}{(1+\|y\|U_1)^{\frac{Q}{n}}})\right)^n
=\left(\E(e^{X/n})\right)^n$$
where $X=\log (1+\|y\|U_1)^{-Q}$.  For fixed $y \in \stackrel{\circ}{B}_d$, $X$ is a bounded random variable. Thus $ \log \left(\E(e^{X/n})\right)^n$ converges to  $\E(X)=-Q\E(\log(1+\|y\|U_1)$, which is Statement 2.

For Statement 3 observe that the sequence $\{W_{n,d}, n=1,2,\ldots\}$ is clearly tight, since all these variables have support in $S_{d-1}$. So take any converging sequence, and call $W^Q_{\infty,d}$ its limit. Since the function $w\mapsto f_y(w)=\frac{1}{(1+\<y,w\>)^Q}$ defined in the unit ball for $y \in \stackrel{\circ}{B}_d$ and $Q>0$ is bounded and continuous, one has that (along this subsequence) $\E\left( \frac{1}{(1+\<y,W_{n,d}\>)^Q}\right)$ converges to $\E\left( \frac{1}{(1+\<y,W^Q_{\infty,d}\>)^Q}\right)$. Therefore
$$
\E\left( \frac{1}{(1+\<y,W^Q_{\infty,d}\>)^Q}\right)=\exp\{QL_d(\| y\|^2\}
$$
and since as mentioned in Section 1 the knowledge of $T_Q(W^Q_{\infty,d})$ characterizes the law of  $W^Q_{\infty,d}$ we have uniqueness of the law of $W^Q_{\infty,d}$ and the stated representation.

For Statement 4 we first observe that Statement 3 says that $W^Q_{\infty,d}$ is the mean vector of a Dirichlet random measure (sometimes called a Dirichlet process) with parameter measure equal to $Q$ times the uniform distribution on $S_{d-1}$. To see this, apply the classical result about the Stieltjes transform of the mean of a Dirichlet random measure due to Cifarelli and Regazzini (1979) and (1990), conveniently stated in Theorem 2.1 of Lijoi and Prunster (2009). Even if this result is stated for a Dirichlet random measure concentrated on the real line, passing from $\R$ to $\R^d$ is standard: it is enough to apply the theorem to the real mean $\<v,W^Q_{\infty,d}\>$ of the Dirichlet process generated by the measure $Q\alpha(du)$ where $\alpha$ is the distribution of $U_1$, $v\in S_{d-1}$ being arbitrary. Finally, recall from Sethuraman (1994) that  $P=\sum_{i=1}^{\infty}\delta_{\Theta_i}\Pi_i$ is actually a Dirichlet random measure governed by $Q$ times the uniform probability on $S_{d-1}.$ Therefore $W_{\infty,d}^Q\sim \int_{S_{d-1}}\theta P(d\theta)$ since both random variables have the same $T_Q$ transform given by Statement 2. Thus the representation  of $W_{\infty,d}^Q$ given in  Statement 4 is obtained. $\square$

The last statement in the previous theorem says that this limit of  Dirichlet walks can be obtained as a walk with an infinite number of steps. This is a particular instance of a more general result apearing in  Hjort and Ongaro (2005) and (2006).

Finally we investigate the random vectors in $\mathbb{R}^d$ whose $Q$-transform has the form $e^{QL_d(\|y\|^2)}$. If $d=1$ we get
$$e^{QL_1(\|y\|^2)}=\frac{1}{(1-\|y\|^2)^{Q/2}}.$$ Applying Proposition 4 , this formula yields that the corresponding limiting distribution has the square radius $(R_{\infty, 1}^{Q})^2$ distributed as $\beta(1/2,Q/2)$ and therefore the distribution of $W^Q_{\infty, 1}=\pm R$ is $$\frac{1}{B(Q/2,Q/2)}(1-w^2)^{\frac{Q}{2}-1}\textbf{1}_{(-1,1)}(w)dw.$$ This result is not easy to get as a limiting case of Proposition 10. 

For $d=2,$ 
formula (\ref{HYP3}) shows that $e^{QL_2(\|y\|^2)}=G(\|y\|^2)^Q$ where $G$ is defined by (\ref{FA}). This time the distribution of $(R_{\infty, 2}^{Q})^2$ is $\beta(1/2,Q)$.

What about $d\geq 3?$  The function $L_d$ when $d$ is odd is explicitly computable. Here is the example $d=3:$ 
$$2L_3(z)=\sum_{k=1}^{\infty}\frac{z^k}{k(1+2k)}=\sum_{k=1}^{\infty}\frac{z^k}{k}-2\sum_{k=1}^{\infty}\frac{z^k}{1+2k}$$
Introducing 
$$f(x)=\sum_{k=1}^{\infty}\frac{x^{2k+1}}{1+2k},\ \ f'(x)=-1+\frac{1}{2(x+1)}-\frac{1}{2(x-1)},\ \ f(x)=-x+\frac{1}{2}\log \frac{1+x}{1-x}$$ we get
\begin{equation}\label{FG3}L_3(\|y\|^2)
=1+\frac{1}{2}\{(1-\|y\|)\log (1-\|y\|)-(1+\|y\|)\log (1+\|y\|)\}\end{equation}

However, finding the distribution of $W^Q_{\infty,3}$ in the unit ball of $\mathbb{R}^3$ such that 
\begin{equation}\label{UFG3}\mathbb{E}\left(\frac{1}{(1+\<y,W^Q_{\infty,3}\>)^Q}\right)=\exp \left(Q L_3(\|y\|^2)\right)\end{equation} or, equivalenty, finding the distribution of  $\|W^Q_{\infty,3}\|^2,$ seems  to be a challenging open problem. Let us give some details about it. Consider the random vector $\Theta$ such that $\Theta$ is  uniform on the unit sphere of $\R^3$ and denote by $U$ its first coordinate. Then $U$  is uniform on $(-1,1)$ by the Archimedes theorem. Denoting $R=\|W^Q_{\infty,3}\|$ we claim that 
$$\mathbb{E}\left(\frac{1}{(1+\<y,W^Q_{\infty,3}\>)^Q}\right)=\E(\left(\frac{1}{(1+\|y\|RU\>)^Q}\right)=\E(f_Q(R\|y\|))$$
where $f_Q(t)=\frac{1}{2}\int_{-1}^1(1+ut)^{-Q}du$ is defined for $t\in (-1,1)$ and  is easy to compute. Replacing for simplicity $\|y\|$ by $t$ in (\ref{FG3}), by using (\ref{UFG3}) the problem of the explicit description of the distribution of $W^Q_{\infty,3}$ is now reduced to the following problem of harmonic analysis: for fixed $Q>0$ find the unique distribution for $R\in [0,1]$ such that for all $t\in (0,1)$
$$\E(f_Q(Rt))=\exp\left(Q\sum_{k=1}^{\infty}\frac{t^{2k}}{2k(2k+1)}\right)=\left(e(1-t)^{\frac{1-t}{2t}}(1+t)^{\frac{-1-t}{2t}}\right)^Q.$$
For instance for $Q=1$ we want $R$ such that
$$\E\left(\frac{1}{2tR}\log \frac{1+tR}{1-tR}\right)=e(1-t)^{\frac{1-t}{2t}}(1+t)^{\frac{-1-t}{2t}}.$$

\section{Dirichlet infinite divisibility and Dirichlet semi groups}
Limits of  Dirichlet random walks or, in view of Theorem 16, means of Dirichlet random measures, are examples of  the following property of infinite divisibility:

\vspace{4mm}\noindent \textbf{Definition.} Let $W\sim \mu$ be a random variable on the unit ball $B_d$ of the Euclidean space of dimension $d.$ We say that $W$ or $\mu$ is \textit{Dirichlet infinitely divisible} of type $Q>0$ if for all $n$ there exists a probability measure $\nu_n$ on the unit ball such that the following occurs: If $Y=(Y_1,\ldots,Y_n)\sim \mathcal{D}(\frac{Q}{n},\ldots,\frac{Q}{n})$ and, independently, $W_1,\ldots,W_n$  are i.i.d. with distribution $\nu_n$, then
$$Y_1W_1+\cdots+Y_nW_n\sim \mu.$$ 
\vspace{4mm}\noindent
Here is the following equivalence property:

\vspace{4mm}\noindent \textbf{Theorem 17.}  The three following properties are equivalent \begin{enumerate}\item $W$ is Dirichlet infinitely divisible of type $Q$;\item for all $n$ there exists a random variable $\tilde W^{(n)}$ on the unit ball such that for all $y\in \stackrel{\circ}{B}_d$ one has 
$$T_Q(W)(y)=\left[T_{Q/n}(\tilde W^{(n)})(y)\right]^n;$$
\item there exists a random variable $\hat W$ on the unit ball such that for all $y\in \stackrel{\circ}{B}_d$ one has  $$T_Q(W)(y)=e^{-Q\mathbb{E}\left(\log (1+\<y,\hat W\>)\right)}.$$
\end{enumerate}

\vspace{4mm}\noindent \textbf{Proof.} For $1)\Rightarrow 2)$ observe that by assumption $W \sim Y_1W^{(n)}_1+\cdots+Y_nW^{(n)}_n$, where we have made explicit the dependence on $n$ of the $W_i$. Thus 
$$
\mathbb{E}\left(\frac{1}{(1+\<y,W\>)^Q}\right)=\mathbb{E}\left(\frac{1
}{(1+\<y,Y_1W^{(n)}_1+\cdots+Y_nW^{(n)}_n\>)^Q}\right)\\
$$
$$
=\mathbb{E}\left[\mathbb{E}\left(\frac{1}{(1+\<y,Y_1W^{(n)}_1+\cdots+Y_nW^{(n)}_n\>)^Q}|W^{(n)}_1,\ldots,W^{(n)}_n\right)\right]\\
$$
$$
=\mathbb{E}\left(\frac{1}{(1+\<y,W^{(n)}_1\>)^{Q/n}}\ldots \frac{1}{(1+\<y,W^{(n)}_n\>)^{Q/n}}\right)
=\left[\mathbb{E}\left(\frac{1}{(1+\<y,W^{(n)}_1\>)^{Q/n}}\right)\right]^n
$$
thus by taking $\tilde W^{(n)}$ equal in law to $W^{(n)}_1$ the representation is obtained.

$2)\Rightarrow 1)$ is obtained since $T_Q(W)$ characterizes the law of $W$. 

$2)\Rightarrow 3)$ Again $\{\tilde W^{(n)}, n=1,2,\ldots\}$ is a tight sequence so we may assume that it has a limit $\hat W$. Arguing as in the previous proof we conclude that $$\lim_{n\rightarrow \infty} \left[\mathbb{E}\left(\frac{1}{(1+\<y,\tilde W^{(n)}\>)^{Q/n}}\right)\right]^n=e^{-Q\mathbb{E}(\log (1+\<y,\hat W\>)}.$$

$3)\Rightarrow 2)$ For any $q>0$, the mean $W^q$ of a Dirichlet random measure with parameter measure $q$ times the law of $\hat W$ has the property
 $$T_q(\hat W^q)(y)=e^{-q\mathbb{E}\left(\log (1+\<y,\hat W\>)\right)}.$$
As a consequence Statement 2 is obtained by taking $\tilde W^{(n)}$ equal in law to $W^{q/n}$.
$\square$

\vspace{4mm}\noindent
Thus from the Dirichlet infinite divisibility property a stronger property follows, via the equivalent Statement 3 of  Theorem 17. Any Dirichlet infinite divisible random variable $W^Q$ of type $Q>0$ is indeed embedded in a Dirichlet semigroup of laws $\{\mu_q, q>0\}$, defined according to the following definition, in the sense that $W^Q \sim \mu_Q$. Notice that this semigroup is weakly continuous at zero since $\mu_q$ converges to the distribution of  $\hat W$ as $q \to 0$.

\vspace{4mm}\noindent \textbf{Definition.} Let $Q\mapsto \mu_Q$ be a map from $(0,\infty)$ to the set of probability measures in the unit ball $B$ of the Euclidean space of dimension $d$. We say that this map is a \textit{Dirichlet semigroup} if for all $n$ and for all $q_1,\ldots,q_n>0$ the following occurs: taking $W_1,\ldots,W_n,Y$ independent such that $W_i\sim \mu_{q_i}$ and $Y=(Y_1,\ldots,Y_n) \sim \mathcal{D}(q_1,\ldots,q_n)$ then setting $Q=q_1+\cdots+q_n$ it holds
$$Y_1W_1+\cdots+Y_nW_n\sim \mu_Q.$$ 
For instance
$$T_Q(W_d^Q)(y)=G(\|y\|^2)^{Q}$$
describes $W_d^Q=R^Q\Theta$, with $\Theta$ uniform in $S_{d-1}$ and $R^2 \sim \beta (\frac{d}{2}, Q+1-\frac{d}{2})$, with $Q>\frac{d}{2}-1$. Thus $W_d^Q$ follows  a hyperuniform law of type $2(Q+1)$. For $d=2$ we get a Dirichlet semigroup and for $d=1$ the restriction to $Q>0$ is a Dirichlet semigroup as well. Similarly the relation

$$
T_Q(\bar W_d^Q)(y)=\frac{1}{(1-\|y\|)^{Q/2}}
$$
describes $\bar W_d^Q=(\bar R)^Q\Theta$, with $\Theta$ uniform in $S_{d-1}$ and $(\bar R)^2 \sim \beta (\frac{d}{2}, \frac {Q+1-d}{2})$, with $Q>d-1$. Thus $\bar W_d^Q$ follows a hyperuniform law of type $Q+1$. Here we get a Dirichlet semigroup only for $d=1$. However by direct inspection of the parameters of these beta distributions it is seen that $\bar W_1^Q \sim W_1^{\frac{Q-1}{2}},$ and with a straightforward calculation one gets $\frac{1+\bar W_1^Q}{2} \sim \beta(\frac{Q}{2},\frac{Q}{2})$. Thus if such a family of symmetric beta distribution is reparametrized by $q=\frac{Q-1}{2}$, for $Q>1$, it keeps the Dirichlet semigroup property.

Notice also that, for any $d>2$ in the first case and any $d>1$ in the second, one still obtains a family of hyperuniform distributions with the Dirichlet semigroup property, but the index of the family runs on a parameter set of the form $(a,+\infty)$, with $a>0$, which clearly remains an additive semigroup.
\section{References}\vspace{4mm} \noindent\textsc{Borwein, J.M., Staub, A., Wan, J.  and  Zudilin, W.} (2012)
'Densities of short uniform random walks.' \textit{Canad. J. Math.} \textbf{64}(5) 961-990
(with an appendix by \textsc{Don Zagier}).

\vspace{4mm} \noindent\textsc{Chamayou, J.-F. and  Letac, G.} (1994)
 'A transient random walk on stochastic matrices with Dirichlet distribution.'
 {\it Ann. Probab.} {\bf
22}: 424-430.

\vspace{4mm} \noindent \textsc{Cifarelli, D.M. and Regazzini, E.} (1979b) 'A general approach to
Bayesian analysis of nonparametric problems. The associative mean values
within the framework of the Dirichlet process'. II. (Italian) \textit{Riv. Mat.
Sci. Econom. Social}. \textbf{2}: 95-111.

\vspace{4mm} \noindent \textsc{Cifarelli, D.M. and Regazzini, E.} (1990) 'Distribution functions of
means of a Dirichlet process'. \textit{Ann. Statist.} \textbf{18}: 429-442 (Correction in
\textit{Ann. Statist.} (1994) \textbf{22}: 1633-1634).

\vspace{4mm}\noindent \textsc{Hjort, N.L. and Ongaro, A.} (2005) 'Exact inference for random Dirichlet means'. \textit{Stat. Inference Stoch. Process.} \textbf{8}: 227-254.

\vspace{4mm}\noindent \textsc{Hjort, N.L. and Ongaro, A.} (2006) 'On the distribution of random Dirichlet jumps'. \textit{Metron} \textbf{64}: 61-92.

\vspace{4mm} \noindent \textsc{Kolesnik, A.D.} (2009) 'The explicit probability distribution of a 6 dimensional random flight' \textit{Theory Stoch. Proc.} \textbf{15}: 33-39.

\vspace{4mm} \noindent \textsc{Le Ca\"{e}r, G. }(2010) 'A Pearson Random Walk with Steps of Uniform Orientation and Dirichlet Distributed Lengths' \textit{J. Stat. Phys.} \textbf{140}: 728-751.

\vspace{4mm} \noindent \textsc{Le Ca\"{e}r, G. }(2011) 'A New Family of Solvable Pearson-Dirichlet Random Walks' \textit{J. Stat. Phys.} \textbf{144}: 23-45.

\vspace{4mm} \noindent\textsc{Lijoi, A. and Pr\"{u}nster, I.} (2009) 'Distributional properties of means of random probability measures' 
\textit{Statist. Surveys} \textbf{3}: 47-95.

 \vspace{4mm} \noindent\textsc{Beghin, L. and Orsingher, E.} (2010) 'Moving randomly amid scattered obstacles' \textit{Stochastics} \textbf{82}: 201-229.

\vspace{4mm} \noindent\textsc{Sethuraman, J.} (1994)
 'A constructive definition of Dirichlet priors.'
 {\it Statist. Sinica} {\bf 4}: 639-650.
\end{document}